A Limit-Free Algebraic–Geometric Construction of the Derivative

with a Foundational Model in the Class of Polynomial Functions

David Kapanadze

Doctor of Pedagogical Sciences

Associate Member, Institute of Mathematics and its Applications (IMA, UK)

Professor, Georgian National University (SEU)

Davitkapanadze06@gmail.com

Abstract

The present study presents an alternative approach to introducing the concept of the derivative in calculus, initially developed in detail within the class of polynomial functions and based on an algebraic–geometric construction without prior use of the formal definition of the limit.

For polynomial functions, the existence and uniqueness of the tangent are characterized by the property that, after subtracting the corresponding linear function from the polynomial, the resulting difference has a double root at the given point. This criterion provides a foundation for a natural construction of the derivative as a new function: to each point corresponds the slope of the tangent constructed at the same point.

Special attention is given to polynomial functions, since in this class it is possible to fully prove algebraically the existence, uniqueness, and basic rules of derivatives. Based on the obtained model, the rules for the derivatives of sums, products, quotients, and composite functions are established, after which the same conceptual scheme is extended to other elementary functions, including exponential, logarithmic, and trigonometric functions.

At a subsequent stage, the constructed model is connected to the linear decomposition of functions, from which the classical limit form of the derivative is naturally obtained. Thus, the limit is presented not as an initial definition but as the analytic expression of the already constructed derivative, creating a conceptual bridge between geometric intuition and formal analysis.

Table of Contents

Introduction

The classical teaching of the derivative concept typically begins with the formal definition of the limit, yet this stage represents one of the first serious conceptual difficulties in instruction. Students are simultaneously required to master both algebraic transformations and abstract limit reasoning, which often leads to the derivative being perceived as a purely formal computational procedure rather than a concept with functional meaning.

To reduce this difficulty, the present study initially selects the class of polynomial functions, since in this class the existence and uniqueness of the tangent can be described clearly with an algebraic criterion. If the corresponding linear function is subtracted from the polynomial graph, the presence of a double root in the resulting difference at a given point expresses the necessary and sufficient condition for tangency.

Based on this criterion, the derivative is introduced as a new function — a correspondence between the points of the argument and the slopes of tangents constructed at the respective points. It is precisely on polynomial functions that a model is created in which the fundamental rules of the derivative can be fully constructed without direct use of the limit. The resulting algebraic–geometric model is subsequently extended to other elementary functions and is connected to the idea of linear decomposition of functions, from which the classical limit formula for the derivative is naturally obtained. This approach demonstrates that the limit can be viewed not as an initial foundation but as the analytic formalization of an already constructed derivative.

1. Pedagogical Overview
 1.1. Pedagogical Foundations

Long-term instructional observations at the upper secondary school level and at the initial stage of higher education have shown that a formal treatment of the topics of the derivative and the limit of a function does not always yield a corresponding conceptual understanding. It was particularly notable that a significant portion of students, despite having previously encountered this material in school, perceived the derivative not as a functional correspondence but rather as a numerical result obtained at an isolated point.
Their knowledge was largely based on formal schemes: the limit was associated with a general notion of infinity, and the derivative with a computational result obtained via the limit. Under such conditions, applying the derivative to practical problems—such as optimization, analysis of functional behavior, or modeling—proved difficult.
These observations indicate that the difficulty is connected not only with the formal definition but also with the initial stage of concept introduction, where the geometric and functional meaning of the derivative has not yet been established for the student.

1.2. Diagnostics and Analysis of Prior Knowledge

To better understand the causes of these difficulties, a systematic observation of the learning process was conducted, analyzing students' prior knowledge, the dynamics of concept acquisition, and the stages where conceptual gaps appeared.
The observation showed that students had firmly mastered several fundamental algebraic and geometric representations: the general notion of a function, the graphical interpretation of linear functions, and the discriminant scheme for solving quadratic equations. Particular significance was attached to the fact that students were well aware of the relationship between the number of roots of a quadratic equation and their graphical location:

- Positive discriminant: two points of intersection;
- Negative discriminant: no intersection;
- Zero discriminant: one double root, which graphically represents the state of tangency.

It was precisely the case of the double root that became the conceptual foundation enabling a natural construction of the tangent in the study of polynomial functions.

1.3. Review of Selected Literature and Differences from Existing Approaches

Attempts to construct the concept of the derivative without limits, or through alternative approaches, are found in the history of mathematics education and analysis in the works of numerous authors (e.g., [1–9]).

Existing studies demonstrate that the directions for alternative constructions of differential calculus are diverse.

Some approaches are based on formal analytic transformations that bypass the definition of the limit (Shisha; Falkowitz, Shisha & Morsi), some are restricted to specific classes of functions, particularly constructions within the class of polynomial functions (Sangwin), while others are oriented toward intuitive or didactic interpretations (Zhang & Tong; Sparks; Tall).

A separate direction is represented by approaches based on infinitesimal quantities (Keisler; Robinson), where the derivative is obtained from the local structure of a function in a particular form.

Another important direction is represented by smooth infinitesimal analysis (Bell), where the derivative is developed within an axiomatic framework based on nilpotent infinitesimals. In this approach, differential calculus is reduced to algebraic manipulation of infinitesimals, and the concept of limit is not required at the initial stage. However, this framework relies on a different foundational system and is conceptually distinct from the algebraic–geometric construction presented in this study.

Conceptually, the approach closest to ours is the representation using dual numbers:

$$f(a + b\varepsilon) = f(a) + f'(a) \cdot b\varepsilon, \qquad \varepsilon^2 = 0,$$

where the derivative appears as the coefficient of the linear term.

The approach presented in this study differs in that the initial construction of the derivative is fully based on the algebraic–geometric properties of polynomial functions. The main foundation is the following fact: the existence and uniqueness of the tangent at a given point is equivalent to the property that, after subtracting the corresponding linear function, the resulting difference has a double root at the same point.

Based on this criterion, the tangent is no longer considered merely as a geometric intuition, but becomes an algebraically constructible object, while the derivative is introduced from the outset as a functional correspondence — a rule that assigns to each point the slope of the tangent.

Unlike existing alternative approaches, the main emphasis here is placed not only on obtaining a specific formula for the derivative, but on constructing an instructional model in which:
• The derivative is perceived from the outset as a function;
• Its existence is justified through an internal algebraic mechanism;
• The fundamental rules of differentiation are derived directly from this construction.

Thus, the class of polynomial functions is presented not as an isolated illustrative case, but as a complete foundational model within which the concept of the derivative can be constructed in a conceptually complete manner without the prior use of the formal definition of the limit.

2. Constructivist Transition: From Old Knowledge to New

Based on prior algebraic knowledge, it is possible to construct a new concept through a natural extension of existing representations.

Consider a polynomial function

$$f(x) = a_n x^n + a_{n-1} x^{n-1} + \cdots + a_1 x + a_0 \qquad (1)$$

and a linear function
$$y = kx + b. \qquad (2)$$

The relative positions of their graphs are determined by the equation

$$f(x) = kx + b. \qquad (3)$$

The number of solutions of this equation determines the number of intersection points between the polynomial graph and the line.

If equation (3) has distinct roots, then the line intersects the polynomial graph at the corresponding points.

A special case arises when one of the roots is repeated.

At this point, a conceptually significant question emerges:
What happens when equation (3) has a double root?

In this case, the line no longer represents merely an intersection line — it touches the polynomial graph at the given point, i.e., it is a tangent.

Thus, for students, a previously known algebraic fact — the existence of a double root — naturally acquires a new geometric interpretation and provides a foundation for constructing the concept of the tangent in the class of polynomial functions.

3. Tangency and the Tangent

3.1. Algebraic Formulation of Tangency

This geometric observation can be formalized algebraically by the following criterion:

If equation (3) has a double root $x = p$, then the line $y = kx + b$ is tangent to the graph of the polynomial function at this point.

Equivalently, the difference
$$f(x) - (kx + b)$$
must be divisible by $(x-p)^2$, i.e.,
$$f(x) - (kx + b) = (x - p)^2 Q(x) \qquad (4)$$

where $Q(x)$ is a polynomial of degree $n - 2$.

Thus, the existence of a double root constitutes a necessary and sufficient algebraic criterion for the existence of a tangent in the class of polynomial functions.

It is precisely this criterion that provides the foundation for the subsequent construction of the derivative, as it allows the slope of the tangent to be determined at each point.

3.2. Existence of the Tangent at Each Point of a Polynomial Function

Consider a polynomial function $f(x)$ and a point on its graph with abscissa $x = p$. The goal is to find a line
$$y = kx + b$$

that is tangent to the graph at this point.

For this, the following condition must be satisfied:

$$f(x) - (kx + b) = (x - p)^2 Q(x) \qquad (5)$$

so that the left-hand side is divisible by $(x-p)^2$.

For every polynomial function, one can always choose parameters $k$ and $b$ such that condition (5) is satisfied algebraically.

Conclusion: At every point of a polynomial function, a tangent exists without invoking the concept of a limit.

3.3. Determination of the Tangent's Slope Coefficient

Once the tangent exists, the question arises: how is its slope determined?

If

$$f(x) - (kx + b) = (x - p)^2 Q(x),$$

then $k$ represents the slope of the tangent at the point $x = p$.

Since $p$ is the common point of the polynomial $f(x)$ and the line $y = kx + b$, it follows that

$$f(p) = kp + b \Rightarrow b = f(p) - kp.$$

Thus, the polynomial admits the local representation

$$f(x) = A(x - p) + f(p) + (x - p)^2 Q(x) \qquad (6)$$

where $A$ is the coefficient of the linear term.

Accordingly,
$$k = A.$$

Therefore, the slope of the tangent is obtained from the local algebraic expansion of the polynomial around the given point.

3.4. Unambiguous Definition of the Tangent

A natural question arises: is it possible to draw a tangent at any point on the polynomial graph?

The answer is affirmative: by an appropriate choice of the parameters $k$ and $b$, condition (5) can always be satisfied.

Fundamental conclusion: At any point on the graph of a polynomial function, a tangent can be drawn, and this tangent is unique.

4. Generation of the Derivative Function

At this stage, a conceptual question arises: do we obtain a new function — a correspondence between the points on the $OX$-axis and the slopes of tangents drawn at the corresponding points?

The answer is affirmative.

Based on the original function, a new function is defined, called the derivative function:

$$f'(x) = k,$$

where $k$ denotes the slope of the tangent at the point $x$.

Concept: The derivative function $f'(x)$ is associated with the original function $f(x)$ such that, at each point, its value equals the slope of the tangent drawn at that point on the graph.

5. Development of the Theory without the Concept of Limit

Based on the algebraic formulation of tangency, we gradually obtain the derivatives of basic functions:

- Derivative of a constant function: $C' = 0$;
- Derivative of a linear function: $(ax + b)' = a$;
- Derivative of a quadratic function: $(x^2)' = 2x$;
- Derivative of a general power function $y = x^n$: $(x^n)' = nx^{n-1}$.

Example 1 — Constant Function

Let $f(x) = C$.

Difference:

$$f(x) - (k(x - x_0) + C) = C - k(x - x_0) - C = -k(x - x_0).$$

This expression is divisible by $(x-x_0)^2$ only when

$$k = 0.$$

Hence,

$$C' = 0.$$

Example 2 — Linear Function

Let

$$f(x) = ax + b.$$

We compute:

$$f(x) - (k(x - x_0) + f(x_0)) = (a - k)(x - x_0).$$

For a double root, it is necessary that
$$a - k = 0 \quad \Rightarrow \quad f'(x) = a.$$

Example 3 — Quadratic Function

Let $f(x) = x^2$.

Difference:

$$x^2 - k(x - x_0) - x_0^2.$$

Rewrite:

$$x^2 - x_0^2 - k(x - x_0) = (x - x_0)(x + x_0 - k).$$

For a double root, we require $x + x_0 - k = 0$ at $x = x_0$, hence:

$$2x_0 - k = 0 \Rightarrow f'(x_0) = 2x_0.$$

Since $x_0$ is arbitrary,

$$f'(x) = (x^2)' = 2x.$$

Example 4 — Cubic Function

Let
$$f(x) = x^3.$$

Difference:
$$x^3 - x_0^3 - k(x - x_0) = (x - x_0)(x^2 + xx_0 + x_0^2 - k).$$

For a double root, we require:
$$3x_0^2 - k = 0 \quad \Rightarrow \quad f'(x_0) = 3x_0^2.$$

Since $x_0$ is arbitrary,
$$f'(x) = (x^3)' = 3x^2.$$

Empirical Conclusion
$$(x^n)' = nx^{n-1}.$$

Fundamental Differentiation Rules

Within the class of polynomials, the basic rules of differentiation are established:

- Sum rule: $(f + g)' = f' + g'$;
- Product rule: $(fg)' = f'g + fg'$;
- Quotient rule: $\left(\frac{f}{g}\right)' = \frac{f'g - fg'}{g^2}$;
- Chain rule: $(f(g(x)))' = f'(g(x)) \cdot g'(x)$.

Proof of the Derivative of a Sum (as an Example)

Let $f$ and $g$ be polynomial functions.

Assume that the tangency condition is satisfied for each function separately.

By adding the corresponding relations, we observe that the tangency condition is also satisfied for the sum, and the slope of the tangent is equal to the sum of the slopes.

Since the point is chosen arbitrarily, it follows that
$$(f + g)' = f' + g'.$$

The theorem is proved.

Remark

The sum of polynomial functions is again a polynomial; however, the objective here is to construct a general model of derivatives.

These results are subsequently extended, from a methodological point of view, to other elementary functions:

- rational and power functions;
- exponential and logarithmic functions;
- trigonometric functions.

Conclusion

Within the class of polynomial functions, the algebraic interpretation of tangency and the construction of the derivative provide a complete and stable theoretical foundation without invoking the concept of limit.

Computer Illustration

An interactive demonstration of the presented method is available at:
https://kapanadze-theory.streamlit.app/

6. Conceptual Transition to the Limit and Classical Analysis

6.1. Formulation of the Problem

At this stage, a natural conceptual question arises:

· Is the constructed theory self-contained within the class of polynomial functions?
· Or does it maintain a connection with classical mathematical analysis, in particular with the concepts of limit and derivative?

These questions are not merely formal; they reveal the necessity of a transition from the algebraic–geometric construction toward the analytic framework of calculus.

Thus, the next conceptual problem may be formulated as follows:
how can the derivative, already constructed as a functional correspondence via tangency, be related to its classical limit representation?

This marks the beginning of a new stage in the development of the theory.

## 6.2. Continuation of Knowledge Construction

To address this problem, the construction of knowledge is continued dialogically with students.

At this stage, the following fundamental concepts are introduced:

- Increment of the argument: $\Delta x$
- Increment of the function: $\Delta y$
- Differential: $dy$

These concepts are not introduced as isolated definitions, but as natural extensions of the previously constructed geometric model.

Their geometric and analytic meanings are then systematically examined, preparing the ground for the transition from the tangent-based interpretation of the derivative to its classical limit form.

## 7. Main Geometric Model

Consider a function
$$y = f(x),$$
and fix a point
$$A(x_0, f(x_0)).$$

Near this point, take a second point
$$B(x_0 + \Delta x, f(x_0 + \Delta x)),$$

where $\Delta x$ represents the increment of the argument.

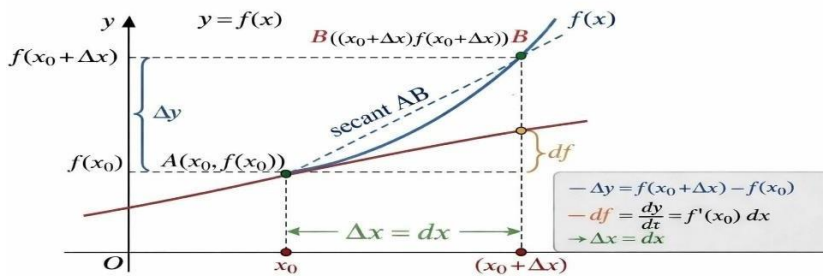

Diagram annotations:

- $A$ — initial point
- $B$ — nearby point on the curve
- $AB$ — secant
- $\Delta x$ — increment of the argument
- $\Delta y$ — increment of the function
- $dy$ — differential

## 8. Definition of Terms

Increment of the argument $\quad \Delta x$
The change of the independent variable.

Increment of the function

$$\Delta y = f(x_0 + \Delta x) - f(x_0)$$

The actual change of the function.

Secant
A secant is a straight line intersecting the curve at two points $A$ and $B$.

Slope of the secant

$$\frac{\Delta y}{\Delta x} = \frac{f(x_0 + \Delta x) - f(x_0)}{\Delta x}.$$

Tangent
As $\Delta x \to 0$, the direction of the secants stabilizes, yielding the tangent at point $A$.

Slope of the tangent

$$f'(x_0).$$

Differential

$$dy = f'(x_0)\,\Delta x.$$

## 9. Linear Decomposition

At this stage, the transition from the geometric interpretation to the analytic formulation is completed.

The total increment of a function can be decomposed as follows:

$$f(x_0 + \Delta x) = f(x_0) + f'(x_0)\, \Delta x + R(\Delta x),$$

where $R(\Delta x)$ is the remainder term.

This decomposition expresses the fundamental idea that the behavior of a function near a point is governed by a linear component (the tangent) together with a higher-order remainder.

## 9.1. Passage to the Limit Form

Dividing the decomposition by $\Delta x$, we obtain:

$$\frac{f(x_0 + \Delta x) - f(x_0)}{\Delta x} = f'(x_0) + \frac{R(\Delta x)}{\Delta x}.$$

Since

$$\frac{R(\Delta x)}{\Delta x} \to 0 \text{ as } \Delta x \to 0,$$

we arrive at the classical limit form:

$$f'(x_0) = \lim_{\Delta x \to 0} \frac{f(x_0 + \Delta x) - f(x_0)}{\Delta x}.$$

## 9.2. Conceptual Remark

Here, the limit does not serve as the initial definition. It appears as:

- a consequence of the existence of the tangent;
- the analytic expression of linear decomposition;
- the vanishing of the higher-order remainder.

Importantly, the derivative does not "approach" the limit — the derivative coincides with this limit.

## 9.3. Conceptual Bridge: From Tangent to Analysis

The geometric model described above shows that the concept of the derivative naturally arises from the tangent problem.

At the initial stage, the change of a function is described by the slope of the secant:

$$\frac{\Delta y}{\Delta x}.$$

As point $B$ approaches $A$, the direction of the secant stabilizes, yielding the tangent with slope

$$f'(x_0).$$

The tangent provides a linear approximation of the function:

$$dy = f'(x_0)\,\Delta x,$$

which describes the local behavior of the function near the point.

The analysis of linear decomposition shows that the actual increment can be written as:

$$f(x_0 + \Delta x) = f(x_0) + f'(x_0)\,\Delta x + R(\Delta x),$$

where

$$R(\Delta x) = o(\Delta x),$$

i.e., the remainder is of higher-order smallness relative to $\Delta x$.

This property leads directly to the classical limit expression:

$$f'(x_0) = \lim_{\Delta x \to 0} \frac{f(x_0 + \Delta x) - f(x_0)}{\Delta x}.$$

Thus, the limit formula is not the initial definition — it emerges as an analytic consequence of the geometric construction of the tangent.

10. General Results and Conclusions

10.1. Pedagogical and Conceptual Outcome

The algebraic–geometric approach to introducing the derivative reveals significant pedagogical and conceptual advantages:

- Derivative as a function from the outset.
  The derivative is perceived not as an isolated numerical value, but as a function assigning to each point the slope of the tangent.
- Introduction without the concept of limit.
  At the initial stage, instruction proceeds without invoking limits, reducing abstraction and facilitating conceptual understanding.
- Preservation of geometric intuition.
  The tangent remains a meaningful geometric object, strengthening both visual and algebraic reasoning.
- Algebraic transparency.
  The slope of the tangent arises directly from local algebraic expansion.
- Conceptual bridge to classical analysis.
  The transition to limits occurs naturally, as an analytic expression of an already constructed concept.

Pedagogical Evaluation

A comparative observation was conducted with two parallel student groups (20 students in each group):

- Control group: taught using the traditional limit-based definition;
- Experimental group: taught using the algebraic–geometric approach.

Observations:

- In the experimental group, the derivative was more often perceived as a functional correspondence;
- Students applied the derivative more flexibly in practical contexts;
- Knowledge was more stable and conceptually organized;
- In the control group, the derivative was often perceived as a formal limit operation.

Conclusion:

The proposed method not only simplifies the initial understanding of the derivative but also promotes a deeper and more functional conceptualization.

10.2. Conceptual Scheme

The approach establishes the following logical sequence:

Graph → Intersection → Tangent → Slope → Functional Correspondence → Derivative

Within this framework, the derivative emerges initially as a function, while the limit appears as a subsequent analytic formalization.

10.3. Main Conclusion

We conclude that:

- The limit is not the foundational starting point of the derivative;
- It is an analytic expression of an already constructed concept;
- The proposed model provides a coherent introduction that preserves:
    - algebraic rigor,
    - geometric intuition,
    - pedagogical consistency.

10.4. Educational Significance

The presented approach:

- reduces initial abstraction;
- strengthens conceptual understanding;
- facilitates transfer to applications;
- remains fully compatible with classical mathematical analysis.